\documentclass[12pt,3p,times]{elsarticle}

\usepackage{srcltx}
\usepackage{xcolor}
\usepackage{colortbl}
\usepackage{enumitem}
\usepackage{amssymb, color}
\usepackage{hyperref}
\usepackage{amsmath}
\usepackage{pifont}
\usepackage{enumitem}
\usepackage[utf8]{inputenc}
\makeatletter
\def\LaTeX{\leavevmode L\raise.42ex
    \hbox{\kern-.3em\size{\sf@size}{0pt}\selectfont A}\kern-.15em\TeX}
\makeatother

\newcommand{\BibTeX}{{\rm B\kern-.05em{\sc
          i\kern-.025emb}\kern-.08em\TeX}}
\journal{}
\makeatletter
\def\@currentlabel{2.1}\label{e:dispaa}
\def\@currentlabel{2.21}\label{e:dispau}
\def\@currentlabel{2.22}\label{e:dispav}
\def\@currentlabel{2.23}\label{e:dispaw}
\def\@currentlabel{2.24}\label{e:dispax}
\def\theequation{\thesection.\@arabic\c@equation}
\makeatother
\renewcommand{\theequation}{\arabic{section}.\arabic{equation}}
\newcommand{\R}{\mathbb R}

\newcommand{\N}{\mathbb N}

\def \O{\Omega}
\everymath{\displaystyle}
\newtheorem{theorem}{Theorem}

\newtheorem{lemma}{Lemma}[section]
\newtheorem{remark}{Remark}[section]
\newtheorem{proposition}{Proposition}[section]

\newtheorem{definition}{Definition}[section]

\renewcommand{\theequation}{\thesection.\arabic{equation}}
\renewcommand{\thesection}{\arabic{section}}
\renewcommand{\theequation}{\thesection.\arabic{equation}}
\let\ssection=\section\renewcommand{\section}{\setcounter{equation}{0}\ssection}
\begin{document}
\begin{frontmatter}
\title{Existence and multiplicity of solutions involving the $p(x)$-Laplacian equations: On the effect of two
nonlocal terms}
\author[mk0,mk1,mk2]{M.K. Hamdani\corref{cor1} }
\ead{hamdanikarim42@gmail.com}
\author[L]{L. Mbarki}
\ead{mbarki.lamine2016@gmail.com}
\author[OMA]{M. Allaoui}
\ead{m.allaoui@uae.ac.ma}
\author[OMA]{O. Darhouche}
\ead{o.darhouche@uae.ac.ma}
\author[dr1,dr2,dr3]{ D.D. Repov\v{s}}
\ead{dusan.repovs@guest.arnes.si}
\begin{center}
\address[mk0] {Science and Technology for Defense Laboratory LR19DN01, CMR, Military Academy, Tunis, Tunisia.}
\address[mk1]{Military Aeronautical Specialities School, Sfax, Tunisia.}
\address[mk2]{Department of Mathematics,  Faculty of Science, University of Sfax, Sfax, Tunisia.}
\address[L]{Department of Mathematics, Faculty of Science, University of Tunis El Manar, Tunis, Tunisia.}
\address[OMA]{Department of Mathematics, FSTH, Abdelmalek Essaadi University, T\'{e}touan, Morocco.}
\address[dr1] {Department of Mathematics, Faculty of Education, University of Ljubljana, Ljubljana, Slovenia}
\address[dr2] {Department of Mathematics, Faculty of Mathematics and Physics, University of Ljubljana, Ljubljana, Slovenia}
\address[dr3] {Department of Mathematics, Institute of Mathematics, Physics and Mechanics,  Ljubljana, Slovenia}
\cortext[cor1]{Corresponding Author: D. D. Repov\v{s}, dusan.repovs@guest.arnes.si}
\end{center}
\begin{abstract}
We study a class of $p(x)$-Kirchhoff
problems
 which is seldom studied because the nonlinearity has nonstandard growth and contains a bi-nonlocal term. Based on variational methods, especially the Mountain pass theorem and Ekeland's variational principle, we obtain the existence of two nontrivial solutions for the problem under certain assumptions. We also apply the Symmetric mountain pass theorem and Clarke's theorem to establish the existence of infinitely many solutions. Our results generalize and extend several existing results.
\end{abstract}

\end{frontmatter}

\section{Introduction}

The purpose  of the present paper is to study the existence  and multiplicity of solutions for the following $p(x)$-Kirchhoff equation, with an additional nonlocal term:
 {\small\begin{eqnarray}\label{Prob1}
\left\{
 \begin{array}{ll}
 M\left(\int_{\Omega}\frac{1}{p(x)}|\nabla u|^{p(x)}dx\right)\Delta_{p(x)}u=\lambda |u|^{p(x)-2}u+ f(x,u)\left[\int_{\Omega}F(x,u)dx\right]^{r}
 \quad \mbox{in }\Omega, \\
 \\
u=0    ~~~~~~~~~~~~~~~~~~~~~~~~~~~~~~~~~~~~~~~~~~~~~~~~~~~~~~~~~~~~~~~~~~~~~~~~~~~~~~~~~~~~~~~~~~~~~~~~~~~ \quad \mbox{on }\partial\Omega,
\end{array}
\right.
\end{eqnarray}}where $\Omega\subset \R^N$ is a bounded smooth domain, $p\in C(\overline{\O}),$
 $N>p(x)>1$, $r>0$ and $\lambda$ are real parameters, $M:\R_0^+\to \R_0^+$ is a Kirchhoff function,  $f:\Omega\times\R\to\R$ is a continuous function satisfying 
 certain
 conditions which will be stated later,
and
 $$F(x,u)=\int_{0}^{u}f(x,t)dt\geq0.$$
 We consider the $p(x)$-Laplacian operator of the form:  $$\Delta_{p(x)}=div(| \nabla u| ^{p(x)-2}\nabla u)=\sum_{i=1}^{N}\left(| \nabla u| ^{p(x)-2}\frac{\partial u}{\partial x_i}\right),
$$
which is not homogeneous and is related to the variable exponent Lebesgue space $L^{p(x)}(\O)$ and the variable exponent Sobolev space $W^{1,p(x)}(\O)$.

These facts imply some difficulties. For example, some classical theories and methods, including the Lagrange multiplier theorem and the theory of Sobolev spaces, cannot be
applied. Problem \eqref{Prob1} is called a bi-nonlocal problem because of the presence of the terms
\[\int_{\Omega}\frac{1}{p(x)}|\nabla u|^{p(x)}dx \mbox{ and } \left[\int_{\Omega}F(x,u)dx\right]^{r},\]
which implies that the first equation in \eqref{Prob1} is no longer a pointwise identity. This phenomenon provokes some mathematical difficulties that make the study of such problems particularly interesting.

Besides, such problems have some physical motivations. Indeed, problem \eqref{Prob1} is related with a physical model introduced by Kirchhoff \cite{K1883} as
follows:
\begin{equation} \label{Eq1.2}
\rho\frac{\partial ^2u}{\partial
	t^2}-\left(a+b \int_0^L\left|\frac{\partial
	u}{\partial x}\right|^2dx\right) \frac{\partial ^2u}{\partial
	x^2}=0,
\end{equation}
where $\rho$, $a$, $b$, $L$ are constants. Here,
 $$M\left(\int_0^L\left|\frac{\partial u}{\partial x}\right|^2dx\right):=a+b\int_0^L\left|\frac{\partial
u}{\partial x}\right|^2dx$$
 describes the
changes of the tension due to the increment in the length of the strings during the
vibrations.

It therefore seems reasonable to be possible to give a realistic meaning for $M(0) = 0$, i.e., when the basic tension of the string is zero. Problem \eqref{Eq1.2} has
 received a lot of attention only after Lions \cite{Lions}  proposed an abstract framework for this problem. We refer the reader to \cite{B1,B2,Chabrowski} for the Laplacian operator and \cite{Choudhuri,CF,H2,J} for the $p$-Laplacian operator.

On the other hand,  there are only a few papers which deal with nonlocal $p(x)$-Kirchhoff equation via variational approach, we can see
\cite{A,AD,AB,CC1,CC2,HHMR,HR,ZGC} and the references therein.
Using variational methods, Corr\^{e}a-Costa \cite{CC1}
investigated the following nonlocal
$p(x)$-Laplacian Dirichlet problem
\begin{eqnarray}\label{Costa}
\begin{cases}
-M\left(\int_{\O}\frac{1}{p(x)}|\nabla u|^{p(x)}dx\right)\Delta_{p(x)} u= h(x,u), \; \mbox{ in }\; \O,\\
u=0 \quad ~~~~~~~~~~~~~~~~~~~~~~~~~~~~~~~~~~~~~~~~~~~~~~~~~~~~~~~ \mbox{on }\partial\Omega,
\end{cases}
 \end{eqnarray}
where $$h(x,u)=\lambda|u|^{q(x)-2}u\left[\int_{\O}\frac{1}{q(x)}|u|^{q(x)}dx\right]^r, \quad
m_0 \leq M(t ) \leq m_1.$$
Here, $m_0$ and $m_1$ are positive constants, and
  $$M(t) = t^{\alpha-1}, \quad
q^-(r+1)<\alpha p^-, \quad
\frac{\alpha(p^+)^\alpha}{(q^-)^{\alpha-1}}<\frac{(q^-)^{r+1}(r+1)}{(q^+)^r}.$$  They proved several results on the existence of positive solutions. Recently, their result 
was extended in Corr\^{e}a-Costa
\cite{CC2} to the general nonlinearities cases: $h(x,u)$ and $M(t)$ were replaced respectively, by
$$f(x,u)\left[\int_{\Omega}F(x,u)\right]^{r}, \quad
Q_1t^{\gamma(x)-1}\leq f(x,t)\leq Q_2 t^{q(x)-1},\quad
A_0+At^{\alpha(x)}\leq M(t)\leq B_0+B t^{\beta(x)},$$  where $A_0, A, B_0, B,Q_1,Q_2$ are positive constants and $\alpha(x),\beta(x),\gamma(x),q(x)\in C_+(\overline{\O})$ satisfy the following conditions
$$\alpha(x)\leq \beta(x)\;\mbox{ and } \gamma(x)\leq q(x)<p^*=\frac{Np(x)}{N-p(x)}.$$
By using Krasnoselskii's genus, they proved the existence of infinitely many solutions for \eqref{Costa}. For a deeper treatment, we refer to \cite{BBE,ZFB} and the references therein.

Motivated by the above results, we are interested in the existence and multiplicity of solutions
for the $p(x)$-bi-nonlocal type problem \eqref{Prob1}. We first state the following conditions for the Kirchhoff
function $M$:
\begin{enumerate}[label=($M\sb{\arabic*}$):, ref=$M\sb{\arabic*}$]
\item\label{M1} $M: [0,+\infty) \rightarrow [0,+\infty)$ is a continuous function such that there exist $t_0\geq0$ and $\gamma\in \left(1,{(p^*)_-}/{p_+}\right)$ satisfying
$$tM(t)\leq \gamma \widehat{M}(t),
\
\hbox{ for all}
\
t\geq t_0,
\mbox{ where }
\
 \widehat{M}(t)=\int_{0}^{t} M(z)dz.$$
 \item\label{M2} There exist positive constants $\alpha, A$ and $C$ such that
$$
\widehat{M}(t)\geq C t^{\alpha}~~\text{for }~t\geq A \geq1~~\text{with}~~\alpha p^->p^+.
$$
\end{enumerate}
A typical prototype of $M$ is given by
\begin{equation}\label{hhj.}
M(t)=a+bt^{\alpha-1},\, \hbox{for all } \ t\geq 0,\,\text{ where } a,b\geq 0,\,b>0 \text{ and } \alpha>1.
\end{equation}
When $M(t) > 0$ for all $t\geq 0$, Kirchhoff problems are said to be nondegenerate and this happens for example if $a>0$ and $b>0$ in the model case \eqref{hhj.}. Otherwise, if $M(0) = 0$ and $M(t) > 0$ for all $t >0$, the Kirchhoff problems are called degenerate and this occurs in the model case \eqref{hhj.} when $a = 0$ and $b>0$.

Moreover, we assume that $f$ is a continuous function which satisfies the following conditions:
\begin{enumerate}[label=($H\sb{\arabic*}$):,
ref=$H\sb{\arabic*}$]
\item\label{H1} The subcritical growth condition holds:
$$
|f(x,s)| \leq C(1+|s|^{q(x)-1}),  \mbox{for all } (x,s) \in \Omega\times \R,
\hbox{where}
\
C>0,
\
p(x)<q(x)<p^*(x); $$
\item\label{H2} The Ambrosetti-Rabinowitz (abbreviated (AR)) condition holds:
$$F(x,s)=\int_{0}^{s}f(x,t)dt$$
is $\theta$-super-homogeneous at infinity, i.e. there exists $s_A > 0$ such that
$$
 0<\theta F(x,s)\leq sf(x,s),\;\mbox{ for all} \; |s| \geq s_A,\; x\in \Omega,
\
\hbox{where}
\
\theta>\frac{\gamma p^+}{r+1};
$$
\item\label{H3} The following holds uniformly in $x \in \O$:
$$ \lim_{s\rightarrow 0} \frac{f(x,s)}{|s|^{p(x)-2}s}=0; $$
\item\label{H4} $f(x,-s)=-f(x,s), \mbox{ for all } (x,s)\in \Omega\times \R$.
\end{enumerate}
\begin{remark}
An example of our conditions being satisfied is given by the following functions:
$M(t)=bt^{\alpha-1}$
and
$f(x,t)=|t|^{q(x)-1}t,$
where
$\alpha>1, b>0$
and
$p(x)<q(x)<p^*(x),$
respectively.
\end{remark}
\begin{remark}
The Ambrosetti-Rabinowitz superlinearity condition was originally introduced by Ambrosetti and Rabinowitz~\cite{AR} and is still used in many works. This condition depicts a superquadratic growth and is used to ensure the boundedness of Palais-Smale sequences of the energy functional and hence in obtaining the mountain
pass geometry. We note that the Palais-Smale condition on the functional is relevant in establishing critical point results and their applications (see also the discussion in \cite{BRW}).
\end{remark}

Now we are in position to state our main results.
\begin{theorem}\label{theorem.1.1}
Suppose that  function $p\in C(\overline{\O})$ satisfies $\gamma p^+<(r+1)p^-$.
Then there exists $\lambda_0>0$ such that for every $\lambda<\lambda_0$, with conditions
\eqref{M1},
\eqref{M2},
\eqref{H1},
\eqref{H2}  and
\eqref{H3}
satisfied,
problem \eqref{Prob1} has at least two nontrivial weak solutions.
\end{theorem}

\begin{theorem}\label{theorem.1.2}
Suppose that  function $p\in C(\overline{\O})$ satisfies $\gamma p^+<(r+1)p^-$.
Then there exists $\lambda_0>0$ such that for every $\lambda<\lambda_0$, with conditions
\eqref{M1},
\eqref{M2},
\eqref{H1},
\eqref{H2},
\eqref{H3}  and
\eqref{H4}
satisfied,
problem \eqref{Prob1} has infinitely many solutions in $W_{0}^{1,p(x)}(\Omega)$.
\end{theorem}

\begin{theorem}\label{theorem.1.3}
Suppose that conditions
 \eqref{M1},
 \eqref{M2},
 \eqref{H1},
 \eqref{H2},
 \eqref{H3}  and
 \eqref{H4}
 are satisfied.
 Then for every $\lambda\in\R$, problem \eqref{Prob1} has infinitely many solutions in $W_{0}^{1,p(x)}(\Omega)$.
\end{theorem}

We conclude with an outline of the structure of the paper. In Section \ref{sec2}, we introduce some
preliminary results concerning  Lebesque and generalized Sobolev spaces and we recall some results that will be used later. In Section \ref{sec3}, we study the Palais-Smale condition. Section \ref{sec4} is devoted to the proof of Theorem \ref{theorem.1.1}. In Section \ref{sec5}, we prove Theorem \ref{theorem.1.2}. Finally, Section \ref{sec6} is dedicated to the proof of Theorem \ref{theorem.1.3}.

\section{Preliminaries}\label{sec2}
In this section, we recall some definitions and basic properties of the generalized Lebesgue space and the variable exponent Sobolev space $W^{1,p(x)}(\Omega)$. For this purpose, let consider $\Omega$ be a bounded domain of $\mathbb{R}^{N}$ and  denote
$$
C_{+}(\overline{\Omega})=\left\{h\in
C(\overline{\Omega})
\mid
h(x)>1, \,\hbox{for all} \ x\in\overline{\Omega}\right\}, 
\quad
h^{+}=\max_{x\in\overline{\Omega}}h(x), \quad
h^{-}=\min_{x\in\overline{\Omega}}h(x), \quad
 h\in C(\overline{\Omega}).
 $$
 The generalized Lebesgue space is defined as
$$L^{p(x)}(\Omega)=\left\{u: u\text{ is a measurable real-valued
function, }\int_{\Omega}|u|^{p(x)}\,dx<\infty\right\}$$
and it is equipped by the following  norm
$$
|u|_{p(x)}=\inf \left\{\lambda>0: \int_{\Omega}
\left|\frac{u(x)}{\lambda}\right|^{p(x)}\,dx\leq 1\right\}.
$$
Thus  $(L^{p(x)}(\Omega),|\cdot|_{p(x)})$ becomes a Banach space.
Let us recall now some results which will be used later.
\begin{proposition}[\cite{Y}] \label{prop2.1}
\begin{itemize}
\item[(1)]  $(L^{p(x)}(\Omega)$, $|\cdot|_{p(x)})$ is  a separable,
uniformly convex Banach space, and has conjugate space
$L^{q(x)}(\Omega)$, where $1/q(x)+1/p(x)=1$.
For every  $u\in L^{p(x)}(\Omega)$ and $v\in L^{q(x)}(\Omega)$, we have
\[
\left|\int_{\Omega}uv\,dx\right|\leq
\left(\frac{1}{p^{-}}+\frac{1}{q^{-}}\right) |u|_{p(x)}|v|_{q(x)}.
\]

\item[(2)] The inclusion between Lebesgue spaces also generalizes the classical framework,
namely, if $0<|\O| <\infty$ and $p_1$, $p_2$ are variable exponents such that $p_1 \leq p_2$ in $\O,$ then there exists a continuous embedding $L^{p_2(x)}(\O)\to L^{p_1(x)}(\O)$.
\end{itemize}
\end{proposition}

An important role in working with the generalized Lebesgue--Sobolev spaces is played by the $m(\cdot)$-modular of the
$L^{p(\cdot)}(\O)$ space, which is the modular $\rho_{p(\cdot)}$ of the space $L^{p(\cdot)}(\O)$
$$
   \rho_{p(\cdot)}(u):=\int_{\Omega} |u|^{p(x)} \,dx.
$$
For more details about these variable exponent Lebesgue spaces see \cite{PRD,RD}.

\begin{lemma}[\cite{FZ2}]\label{s+}
Denote
$$
\Lambda(u)=\int_\O\frac{1}{p(x)}|\nabla u| ^{p(x)}dx,\;\mbox{ for all } u\in W_0^{1,p(x)}(\O).
$$
Then $\Lambda(u)\in C^1(W_0^{1,p(x)}(\O),\R)$  and the derivative operator $\Lambda'$ of $\Lambda$ is
$$
\langle \Lambda'(u),v\rangle=\int_\O|\nabla u|^{p(x)-2}\nabla u \nabla vdx,\;\mbox{ for all } u,v\in W_0^{1,p(x)}(\O),
$$
and the following holds:
\begin{enumerate}
  \item $\Lambda$ is a convex functional;
  \item $\Lambda': W_0^{1,p(x)}(\O)\to (W^{-1,p'(x)}(\O))=\left(W_0^{1,p(x)}(\O)\right)^*$ is a bounded homeomorphism and strictly monotone operator, and the conjugate exponent satisfies $\frac{1}{p(x)}+\frac{1}{p'(x)}=1$;
  \item $\Lambda'$ is a mapping of type $(S_+)$, namely, $u_n\rightharpoonup u$ and $\limsup \langle \Lambda'(u_n),u_n-u\rangle\leq 0$,
imply $u_n\to u$ (strongly) in $W_0^{1,p(x)}(\O)$.
\end{enumerate}
\end{lemma}
\begin{definition}
 We say that $u\in W_0^{1,p(x)}(\O)$ is a weak solution of problem \eqref{Prob1}, if
\[
M\left(\int_{\Omega}\frac{1}{p(x)}|\nabla u|^{p(x)}dx\right)\int_{\Omega}|\nabla u|^{p(x)-2}\nabla u \nabla v dx-\lambda\int_{\Omega}|u|^{p(x)-2}uv dx=\]
\[
 \left[\int_{\Omega}F(x,u)dx\right]^{r}\int_{\Omega}f(x,u)vdx,
\
\hbox{where}
\
 v\in W_0^{1,p(x)}(\O).
 \]
\end{definition}
The energy functional $J_\lambda: W_0^{1,p(x)}(\O) \to \mathbb{R}$ associated with problem \eqref{Prob1}
\begin{eqnarray}\label{funct energ}
J_\lambda(u)&=&
\widehat{M}\int_{\Omega}
\frac{|\nabla u|^{p(x)}}{p(x)}dx
-\lambda\int_\Omega\frac{| u|^{p(x)}}{p(x)}dx-\frac{1}{r+1}\left[\int_{\Omega}F(x,u)dx\right]^{r+1}\nonumber\\
&:=&\Phi(u)-E_\lambda(u)-\Psi(u),
\
\hbox{for all}
\
u\in W_0^{1,p(x)}(\O),
\end{eqnarray}
is well-defined and of $C^1$-class on $W_0^{1,p(x)}(\O)$. Moreover, we have
\begin{eqnarray}\label{deriv funct energ}
\langle J_\lambda'(u),v\rangle&=&M\int_{\Omega}\frac{|\nabla u|^{p(x)}}{p(x)}dx
\int_{\Omega}|\nabla u|^{p(x)-2}\nabla u \nabla v dx-\lambda\int_{\Omega}|u|^{p(x)-2}uv dx\nonumber\\&&-\left[\int_{\Omega}F(x,u)dx\right]^{r}\int_{\Omega}f(x,u)vdx,
\
\hbox{for all}
\
u, v \in W_0^{1,p(x)}(\O).
\end{eqnarray}
Hence, we can observe that the critical points of  functional $J_\lambda$ are the weak
solutions for problem \eqref{Prob1}. In order to simplify the presentation we will denote the norm of $W_0^{1,p(x)}(\O)$ by $\| .\| $ instead of $\| \cdot\| _{W_0^{1,p(x)}(\O)}$.
For simplicity, we use $C_i, i=1,2,...$ to denote general positive constants whose exact values may change from one place to another.
\section{The Palais-Smale compactness condition}\label{sec3}
\begin{definition}
\label{def1} Let $(W_0^{1,p(x)}(\O),\;\|.\| )$ be a Banach space and $J_\lambda \in C^1(W_0^{1,p(x)}(\O))$. Given $c \in\R$, we say that $J_\lambda$ satisfies the Palais--Smale condition at the level $c\in \R$ (``$(PS)_c$ condition'' for short) if any sequence $(u_n) \in W_0^{1,p(x)}(\O)$ satisfying
\begin{equation}
\label{condps}
 J_\lambda(u_n) \rightarrow c \mbox{ and } J_\lambda'(u_n)\rightarrow  0  \mbox{ in } W^{-1,p'(x)}(\O) \mbox{ as } \;n\rightarrow \infty,
\end{equation}
has a convergent subsequence.
\end{definition}

\begin{lemma}\label{PS}
Assume that  conditions
\eqref{M1},\eqref{M2},\eqref{H1}
and \eqref{H2} are satisfied. Then  functional  $J_\lambda$ satisfies the $(PS)_c$ condition for any $c\neq0$.
\end{lemma}
{\bf  Proof.}  We proceed in two steps.

\vspace*{4pt}\noindent
\textbf{Step 1.} We prove that $(u_n)$ is bounded in $W_0^{1,p(x)}(\O)$. Let $(u_n) \subset W_0^{1,p(x)}(\O)$ be a $(PS)_c$ sequence for any $c\neq0$.
By \eqref{M1}, for $\|u\|$ large enough,
\begin{equation}\label{al1.1}
\begin{split}
  \gamma p^+ \Phi(u)  =\gamma p^+\widehat{M}(\Lambda(u))\geq p^+M(\Lambda(u))\Lambda(u) \geq M(\Lambda(u))\int_\Omega|\nabla u|^{p(x)}dx=\Phi'(u)u.
\end{split}
  \end{equation}
From \eqref{H2} we can see that there exists $C_1>0$ such that
$$
-C_1\leq \theta \int_{\Omega}F(x,u)dx\leq \int_{\Omega}f(x,u)udx+C_1,~~\mbox{for all }u\in W_0^{1,p(x)}(\O),
$$
and thus, given any $\varepsilon\in(0,\theta)$, there exists $A_\varepsilon\geq A$ such that
\begin{equation}\label{al1.2}
  (\theta-\varepsilon)\int_{\Omega}F(x,u)dx\leq \int_{\Omega}f(x,u)udx~~\text{if} ~\int_{\Omega}F(x,u)dx\geq A_\varepsilon.
\end{equation}
We may assume $A_\varepsilon>\frac{C_1}{\theta}$. Note that in this case, the inequality $\int_{\Omega}F(x,u)dx\geq A_\varepsilon$ is equivalent to $\left|\int_{\Omega}F(x,u)dx\right|\geq A_\varepsilon,$ because
$\int_{\Omega}F(x,u)dx\geq-\frac{C_1}{\theta},$
 for all
$u\in W_0^{1,p(x)}(\O).$
We claim that there exists $C_\varepsilon>0$ such that
\begin{equation}\label{al1.3}
  \Psi'(u)u-(r+1)(\theta-\varepsilon)\Psi(u)\geq-C_\varepsilon,~~\hbox{for all } u\in W_0^{1,p(x)}(\O).
\end{equation}
Indeed, if $\left|\int_{\Omega}F(x,u)dx\right|\leq A_\varepsilon$, then the validity of \eqref{al1.3} is obvious. When
$\left|\int_{\Omega}F(x,u)dx\right|\geq A_\varepsilon,$
i.e.,
$\int_{\Omega}F(x,u)dx\geq A_\varepsilon,$
 it follows by \eqref{al1.2}
that
\begin{align*}
  (r+1)(\theta-\varepsilon)\Psi(u) & =(\theta-\varepsilon)\left(\int_{\Omega}F(x,u)dx\right)^{r+1}
=(\theta-\varepsilon)\left(\int_{\Omega}F(x,u)dx\right)^{r}\int_{\Omega}F(x,u)dx\\
&   \leq \left(\int_{\Omega}F(x,u)dx\right)^{r} \int_{\Omega}f(x,u)udx
=\Psi'(u)u,
\end{align*}
and so \eqref{al1.3} holds.

Now let $(u_n)\subset W_0^{1,p(x)}(\O)\backslash\{0\}$, $J_\lambda'(u_n)\to 0$ and $J_\lambda(u_n)\to c$ with $c\neq0$.
Since $\gamma p^+<(r+1)\theta$, there exists $\varepsilon>0$ small enough
so  that $\gamma p^+<(r+1)(\theta-\varepsilon)$.
Then, since $(u_n)$ is a $(PS)_c$  sequence, applying \eqref{al1.1}, \eqref{al1.3} and \eqref{M2}, for sufficiently large $n$ we have
$$
    (r+1)(\theta-\varepsilon)c+1+\|u_n\|
    \geq (r+1)(\theta-\varepsilon)J_\lambda(u_n)-J_\lambda'(u_n)u_n
$$
$$
 \geq \left((r+1)(\theta-\varepsilon)-\gamma p^+\right)\Phi(u_n)+\left(\gamma p^+\Phi(u_n)-\Phi'(u_n)u_n \right)
 $$
$$
+\left(\Psi'(u_n)u_n-(r+1)(\theta-\varepsilon)\Psi(u_n)\right)
-\lambda(r+1)(\theta-\varepsilon)\int_\Omega\frac{1}{p(x)}|u_n|^{p(x)}dx
$$
$$
+\lambda\int_\Omega|u_n|^{p(x)}dx
\geq C_2\|u_n\|^{\alpha p^-}-C_3-C_\varepsilon-
\lambda\int_\Omega\left(\frac{(r+1)(\theta-\varepsilon)}{p(x)}-1\right)|u_n|^{p(x)}dx
$$
$$
\geq \left\{
        \begin{array}{ll}
          C_2\|u_n\|^{\alpha p^-}-C_3-C_\varepsilon, & \hbox{if $\lambda\leq 0$} \\
          C_2\|u_n\|^{\alpha p^-}-C_3-C_\varepsilon-\lambda\left(\frac{(r+1)(\theta-\varepsilon)}{p^-}-1\right)C_4\|u_n\|^{p^+}, & \hbox{if $\lambda>0$.}
        \end{array}
      \right.
      $$
Since $\alpha p^->p^+>1$,
the above inequalities 
imply that $(u_n)$ is
bounded in
 $W_0^{1,p(x)}(\O)$.

\vspace*{4pt}\noindent\textbf{Step 2.} Now we claim that $(u_n)$ has a strongly convergent subsequence. To complete the argument we need the following proposition.
\begin{proposition}
\begin{enumerate}
  \item[$(i)$] Functional $\Phi:X:=W_0^{1,p(x)}(\O) \to\mathbb{R} $ is sequentially weakly lower semi-continuous, $\Psi, E_\lambda:X \to\mathbb{R}$ are sequentially weakly continuous, and
thus $J_\lambda$ is sequentially weakly lower semi-continuous.
  \item[$(ii)$] Mappings $\Psi',E_\lambda':X \to X^*$ are sequentially weakly-strongly continuous. For any open set $D\subset X\backslash\{0\}$ with $\overline{D}\subset X\backslash\{0\}$, mappings $\Phi'$ and $J_{\lambda}':\overline{D}\to X^*$ are bounded and  of type $(S+)$.
\end{enumerate}
\end{proposition}
{\bf  Proof.} 
\begin{enumerate}
  \item [${(i)}$] Since  function $\widehat{M}(t)$ is increasing and  functional $\Lambda$ is
sequentially weakly lower semi-continuous, we can see that  functional $\Phi:X:=W_0^{1,p(x)}(\O) \to\mathbb{R} $ is sequentially weakly lower
semi-continuous.
\item [${(ii)}$]  Noting that  embedding $X\hookrightarrow L^{q(x)}(\Omega)$  is compact, we can see that $\Psi$, $\Psi'$, $E_\lambda,$ and $E_\lambda'$ are sequentially
weakly-strongly continuous. Now let $\overline{D}\subset X\backslash\{0\}.$ It is clear that mappings $\Phi'$ and $J_\lambda':\overline{D}\to X^*:=(W^{-1,p'(x)}(\O))$ are bounded.
 To prove that $\Phi':\overline{D}\to X^*$ is of type $(S+)$, assume that $(u_n)\subset \overline{D}$, $u_n\rightharpoonup u$ in $X$ and $\limsup_{n\to+\infty}\Phi'(u_n)(u_n-u)\leq0$.
 Then there exist positive constants $C_1$ and $C_2$ such that
$C_1\leq \Lambda(u_n) \leq C_2$ and therefore there exist positive constants $C_3$ and $C_4$ such that $C_3\leq M(\Lambda(u_n)) \leq C_4$. Noting that $\Phi'(u_n)=M(\Lambda(u_n))\Lambda'(u_n)$, it follows from $\limsup_{n\to+\infty}\Phi'(u_n)(u_n-u)\leq0$ that $\limsup_{n\to+\infty}\Lambda'(u_n)(u_n-u)\leq0$.
Since $\Lambda'$ is of type $(S+)$, we obtain $u_n\to u$  in $X$. This shows that  mapping $\Phi':\overline{D}\to X^*$ is
of type $(S+)$. Moreover, since $\Psi'$ and $E_\lambda'$ are sequentially weakly-strongly continuous,  mapping $J_{\lambda}':\overline{D}\to X^*$ is of type $(S+)$.\qed   
\end{enumerate}

We can now complete the proof of Step 2. Since $J_{\lambda}(0)=0$ and $J_{\lambda}(u_n)\to  c\neq0$, there exists $\varepsilon>0$ small enough such that for sufficiently large $n$, $\|u_n\|>\varepsilon$.
Setting $D=\{u\in W_{0}^{1,p(x)}(\Omega)\,/\,\|u_n\|>\varepsilon\}$, then $u_n\in D$ for $n$ sufficiently large. Because $(u_n)$ is bounded, we can consider a subsequence of $(u_n)$, still denoted by $(u_n)$, such that $u_n\in D$ and $u_n\rightharpoonup u$. The condition $J_{\lambda}'(u_n)\to 0$ implies $J_{\lambda}'(u_n)(u_n-u)\to 0$. Since $J_{\lambda}': \overline{D}\to  W_{0}^{1,p(x)}(\Omega)^{*}$ is of $(S+)$ type, we have $u_n\to  u\in\overline{D}$.
\qed   
\section{Proof of Theorem \ref{theorem.1.1}}\label{sec4}
\label{section4}
In this section, the existence of nontrivial weak solutions for \eqref{Prob1} is shown by applying
the Mountain pass theorem and a variant of the Ekeland variational principle under suitable
assumptions.
To verify the conditions of the Mountain pass theorem (see e.g.,~\cite{Willem}),
we first need to prove two lemmas.
\begin{lemma}\label{negatif}
Suppose that conditions
 \eqref{M1},\eqref{H1}
 and \eqref{H2} are satisfied. Then for any $w\in W_{0}^{1,p(x)}(\Omega)\backslash\{0\}$,
$J_{\lambda}(sw)\rightarrow-\infty \mbox{ as } s\rightarrow+\infty.$
\end{lemma}
{\bf  Proof.}  Let $w\in W_{0}^{1,p(x)}(\Omega)\backslash\{0\}$ be given. From \eqref{M1} and for $t\geq1$, we can easily obtain that
 $\widehat{M}(t)\leq \widehat{M}(1)t^{\gamma}.$
Then $$E(sw)=\widehat{M}\left(\int_{\Omega}\frac{1}{p(x)}|\nabla sw|^{p(x)}dx\right)\leq d_{1}s^{\gamma p^+},$$
for $s$ large enough and  $d_{1}$  a positive constant depending on $w$.
By conditions \eqref{H1} and \eqref{H2}, we have
$$\left[\int_{\Omega}F(x,sw)dx\right]^{r+1}\geq d_{2}s^{(r+1)\theta},$$ for $s$ large enough and where $d_{2}$ is a positive constant depending on $w$.
Finally, we have
\begin{equation*}
  \left|\int_\Omega\frac{ 1}{p(x)}| sw|^{p(x)}dx\right|\leq\frac{1}{p^+}\left(\int_{\Omega}{|w|^{p(x)}dx}\right)s^{p+}=d_{3}s^{p^+},
\end{equation*}
for $s$ large enough, where $d_{3}$ is a positive constant depending on $w$.
Hence for any $w\in W_{0}^{1,p(x)}(\Omega)\backslash\{0\}$ and $s$ large enough,
\begin{equation*}
J_{\lambda}(sw)\leq\left\{\begin{array}{lll}
d_{1}s^{\gamma p^+}-d_{2}s^{(r+1)\theta}+\lambda d_{3}s^{p^+} &\mbox{if} &\lambda>0, \\
d_{1}s^{\gamma p^+}-d_{2}s^{(r+1)\theta}-\lambda d_{3}s^{p^+} &\mbox{if} &\lambda\leq0.
\end{array}\right.
\end{equation*}
Thus, since $p^{+}\leq\gamma p^{+}<(r+1)\theta$, we conclude that $J_{\lambda}(sw)\to-\infty$ as $s\to+\infty$.
\qed   
\begin{lemma}\label{geometryMPT}
Suppose that conditions  \eqref{M1},\eqref{H1}
and \eqref{H3} are satisfied. Then there exist positive numbers $a,\rho,\lambda_0$ such that
$
J_\lambda(u)\geq a>0~~\text{if}~\|u\|=\rho ~~\text{and}~ \lambda<\lambda_0.
$
\end{lemma}
{\bf  Proof.}  Conditions \eqref{H1} and \eqref{H3} imply that
$
|F(x,t)|\leq \varepsilon|t|^{p(x)}+C_{\varepsilon}|t|^{q(x)}, ~\text{ for all}~(x,t)\in\Omega\times\mathbb{R}.
$
For $\|u\|$ small enough, we have
\begin{align*}
  \int_\Omega F(x,u)dx & \leq \varepsilon\int_\Omega|u|^{p(x)}dx+C_{\varepsilon}\int_\Omega|u|^{q(x)} dx
    \leq \varepsilon\left(|u|_{p(x)}^{p^+}+|u|_{p(x)}^{p^-}\right)
    +C_{\varepsilon}\left(|u|_{q(x)}^{q^+}+|u|_{q(x)}^{q^-}\right)\\
    & \leq \varepsilon\left(C_1^{p^+}\|u\|^{p^+}+C_1^{p^-}\|u\|^{p^-}\right)+C_{\varepsilon}\left(C_2^{q^+}\|u\|^{q^+}+C_2^{q^-}\|u\|^{q^-}\right)\\
    & \leq \varepsilon\left(C_1^{p^+}+C_1^{p^-}\right)\|u\|^{p^-}+C_{\varepsilon}\left(C_2^{q^+}+C_2^{q^-}\right)\|u\|^{q^-}\\
     & \leq \varepsilon\left(C_1^{p^+}+C_1^{p^-}\right)\|u\|^{p^-}+C_{\varepsilon}\left(C_2^{q^+}+C_2^{q^-}\right)\|u\|^{p^-}
      \leq C_3\|u\|^{p^-},
\end{align*}
where $C_3=\varepsilon\left(C_1^{p^+}+C_1^{p^-}\right)+C_{\varepsilon}\left(C_2^{q^+}+C_2^{q^-}\right).$
Therefore
\begin{equation}\label{al3.1}
  \Psi(u)=\frac{1}{r+1}\left[\int_{\Omega}F(x,u)dx\right]^{r+1}\leq\frac{C_3^{r+1}}{r+1}\|u\|^{(r+1)p^-}.
\end{equation}
Moreover, condition \eqref{M1} gives
\begin{equation}\label{al3.2}
  \widehat{M}(t)\geq\widehat{M}(1)t^\gamma, ~~\text{for all}~t\in [0,1].
\end{equation}
Thus, using \eqref{al3.1} and \eqref{al3.2}, we obtain
\begin{equation}\label{al3.3}
\begin{split}
  J_\lambda(u) & = \widehat{M}(\Lambda(u))-\lambda\int_\Omega\frac{ 1}{p(x)}| u|^{p(x)}dx-\Psi(u)\\
   & \geq \left\{
            \begin{array}{ll}
              \widehat{M}(1)\left(\Lambda(u)\right)^\gamma-\frac{C_3^{r+1}}{r+1}\|u\|^{(r+1)p^-}, & \hbox{if $\lambda\leq0$} \\
              \widehat{M}(1)\left(\Lambda(u)\right)^\gamma-\frac{\lambda C_1^{p^-}}{p^-}\|u\|^{p^-}-\frac{C_3^{r+1}}{r+1}\|u\|^{(r+1)p^-}, & \hbox{if $\lambda>0$}
            \end{array}
          \right.\\
& \geq \left\{
            \begin{array}{ll}
              \frac{\widehat{M}(1)}{\left(p^+\right)^{\gamma}}\|u\|^{\gamma p^+}-\frac{C_3^{r+1}}{r+1}\|u\|^{(r+1)p^-}, & \hbox{if $\lambda\leq0$} \\
              \frac{\widehat{M}(1)}{\left(p^+\right)^{\gamma}}\|u\|^{\gamma p^+}-\frac{\lambda C_1^{p^-}}{p^-}\|u\|^{p^-}-\frac{C_3^{r+1}}{r+1}\|u\|^{(r+1)p^-}, & \hbox{if $\lambda>0$}
            \end{array}
          \right.\\
& = \left\{
            \begin{array}{ll}
             \|u\|^{\gamma p^+}\left( \frac{\widehat{M}(1)}{\left(p^+\right)^{\gamma}}-\frac{C_3^{r+1}}{r+1}\|u\|^{(r+1)p^--\gamma p^+}\right), & \hbox{if $\lambda\leq0$} \\
             \|u\|^{\gamma p^+}\left( \frac{\widehat{M}(1)}{\left(p^+\right)^{\gamma}}-\frac{\lambda C_1^{p^-}}{p^-}\|u\|^{p^--\gamma p^+}-\frac{C_3^{r+1}}{r+1}\|u\|^{(r+1)p^--\gamma p^+}\right), & \hbox{if $\lambda>0$.}
            \end{array}
          \right.\\
\end{split}
\end{equation}
Now, for each $\lambda>0$, we define a continuous function $h_\lambda : (0,\infty)\to\mathbb{R}$,
$$
h_\lambda(t)=\frac{\lambda C_1^{p^-}}{p^-}t^{p^--\gamma p^+}+\frac{C_3^{r+1}}{r+1}t^{(r+1)p^--\gamma p^+}.
$$
Since $1<p^-<\gamma p^+<(r+1)p^-$, it follows that
$
\lim\limits_{t\to0^+}h_\lambda(t)=\lim\limits_{t\to +\infty}h_\lambda(t)=+\infty.
$
Thus we can find the infimum of $h_\lambda(t)$. Note that equating
$$
h'_\lambda(t)=\frac{\lambda C_1^{p^-}(p^--\gamma p^+)}{p^-}t^{p^--\gamma p^+-1}+\frac{C_3^{r+1}((r+1)p^--\gamma p^+)}{r+1}t^{(r+1)p^--\gamma p^+-1}=0,
$$
we get $$t_0=t=C_4\lambda^{\frac{1}{rp^-}},
\
\hbox{where}
\
C_4=\left(\frac{c_1^{p^-}(r+1)(\gamma p^+-p^-)}{C_3^{r+1}p^-((r+1)p^--\gamma p^+)}\right)^{\frac{1}{rp^-}}>0.
$$
Clearly, $t_0>0$. It can also be checked that $h_\lambda''(t_0)>0$
 and hence the infimum of $h_\lambda(t)$ is achieved at $t_0$. Now, observing that
$$
h_\lambda(t_0)=\left(\frac{C_1^{p^-}C_4^{p^--\gamma p^+}}{p^-}+\frac{C_3^{r+1}C_4^{(r+1)p^--\gamma p^+}}{r+1}\right)\lambda^{\frac{(r+1)p^--\gamma p^+}{rp^-}}\to 0 ~~\text{as }~\lambda\to0^+,
$$
we can infer from \eqref{al3.3} that there exists $\lambda_0>0$ such that  for all $\lambda<\lambda_0$ we can choose
$\rho$ small enough and $a>0$ such that
$
J_\lambda(u)\geq a>0, ~~\text{for all }~u\in X  ~~\text{with}~\|u\|=\rho.
$
\qed   

Let $\lambda_0>0$ be a constant as given in Lemma \ref{geometryMPT}.
By Lemmas \ref{PS}, \ref{negatif}, \ref{geometryMPT} and the Mountain pass theorem, we deduce that for all  $\lambda\in (0,\lambda_0)$, $J_\lambda$ has a critical point $u_{1}\in X$ which is a weak solution for problem \eqref{Prob1}. Moreover, $u_{1}$ satisfies
\begin{equation}\label{positif}
J_\lambda(u_{1})\geq a>0,
\end{equation}
which implies that $u_{1}$ is nontrivial.

 We will show that there exists a second weak solution $u_2 \ne u_1$ by using the Ekeland variational principle. By Lemma \ref{geometryMPT}, we have
$
\inf_{u\in \partial B(0,r)}(J_\lambda(u) )>0,
$
and
by Lemma \ref{negatif}, there exists $w\in X$ such that $J_\lambda(tw)<0$ for $t>0$ large enough. Moreover, as in the proof of Lemma \ref{geometryMPT}, for $u\in B(0,r),$ we have
$$
J_{\lambda}(u)\geq  \left\{
            \begin{array}{ll}
             \|u\|^{\gamma p^+}\left( \frac{\widehat{M}(1)}{\left(p^+\right)^{\gamma}}-\frac{C_3^{r+1}}{r+1}\|u\|^{(r+1)p^--\gamma p^+}\right),  ~~\hbox{if} ~~\lambda\leq 0 \\
             \|u\|^{\gamma p^+}\left( \frac{\widehat{M}(1)}{\left(p^+\right)^{\gamma}}-\frac{\lambda C_1^{p^-}}{p^-}\|u\|^{p^--\gamma p^+}-\frac{C_3^{r+1}}{r+1}\|u\|^{(r+1)p^--\gamma p^+}\right),~~\hbox{ if}~~ \lambda>0.
            \end{array}
          \right.
$$
Therefore
$
-\infty <\underline{c}=\inf_{u\in \overline{B(0,r)}}(J_\lambda(u) )<0.
$

Let $\varepsilon>0$, be such that
$0<\varepsilon <\inf_{u\in \partial B(0,r)}(J_\lambda(u) )-\inf_{u\in B(0,r)}(J_\lambda(u) ).$
We deduce from the above information that  functional $J_\lambda :\overline{B(0,r)}\rightarrow \mathbb{R}$, is lower bounded and $J_\lambda \in C^{1}(\overline{B(0,r)},\mathbb{R})$. Therefore, by using the Ekeland principle, we conclude that there exists  $u_{\varepsilon }\in
\overline{B(0,r)}$, such that
\begin{equation*}
\left\{
\begin{array}{ll}
\underline{c}\leq J_\lambda (u_{\varepsilon })\leq \underline{c}+\varepsilon  \\
\\
J_\lambda (u_{\varepsilon })<J_\lambda (u)+\varepsilon ||u-u_{\varepsilon }||,u\neq
u_{\varepsilon } .
\end{array}%
\right.
\end{equation*}%
Since
$
J_\lambda (u_{\varepsilon })\leq \inf_{u\in \overline{B(0,r)}}(J_\lambda(u))+\varepsilon \leq\inf_{B(0,r)}(J_\lambda(u) )+\varepsilon <\inf_{\partial B(0,r)}(J_\lambda(u)),
$
we can  deduce that $u_{\varepsilon }\in B(0,r)$.

Now, we define $\Xi_\lambda:\overline{B(0,r)}\rightarrow \mathbb{R}$
 by
 $\Xi_\lambda(u)=J_\lambda (u)+\varepsilon \|u-u_{\varepsilon }\|.$
It is clear that $u_{\varepsilon }$ is a minimum of $\Xi_\lambda$. Therefore,  for $t>0$ large enough and for any $v\in B(0,1),$ we have
$$ \frac{\Xi_\lambda(u_{\varepsilon }+tv)-\Xi_\lambda(u_{\varepsilon})}{t}\geq 0,
\
\hbox{that is,}
\
 \frac{J_\lambda (u_{\varepsilon }+tv)-J_\lambda (u_{\varepsilon })}{t}+\varepsilon\left\Vert v\right\Vert \geq 0.$$
By letting $t$ tend to infinity, we obtain
$
J'_\lambda(u_{\varepsilon })(v)+\varepsilon \left\Vert v\right\Vert
\geq 0.
$
This implies that
$ 
\left\Vert J'_\lambda(u_{\varepsilon })\right\Vert \leq \varepsilon .
$ 
By the argument above, we deduce the existence of  a sequence $
(u_n)\subset B(0,r)$, such that
\begin{equation}\label{estimation1}
J_\lambda (u_{n})\rightarrow \underline{c}<0,\;\text{ and }\;J'_\lambda(u_{n})\rightarrow 0.
\end{equation}%
Since  $
 (u_{n}) \subset B(0,r)$, it follows that $(u_{n})$ is bounded in $X$. So, up to a subsequence, there exists $u_2\in X$ such that  $(u_{n}) $
converges weakly to $u_2\in X$. Hence, by the proof of Lemma \ref{PS}, we deduce that  $u_{n}\rightarrow u$ strongly in $X.$
\\Since $J_\lambda \in C^{1}(X,\mathbb{R})$, we have
$ 
J'_\lambda(u_{n})\rightarrow J'_\lambda(u_2),\text{ as }%
n\rightarrow \infty.
$ 
Hence, from \eqref{estimation1}, we conclude that
\begin{equation}\label{estimation2}
J'_\lambda(u_2)=0,\;\left\Vert u_2\right\Vert <r,\text{ and }%
J_\lambda (u_2)<0.
\end{equation}%
This implies that $u_2$ is a nontrivial solution for problem \eqref{Prob1}.
Finally, by combining \eqref{positif} and \eqref{estimation2}, we obtain $J_\lambda (u_2)<0<J_\lambda (u_1).$
The proof of Theorem \ref{theorem.1.1} is now completed.\qed
\section{Proof of Theorem \ref{theorem.1.2}}\label{sec5}
In this section, we will show that 
problem \eqref{Prob1} has infinitely many pairs $(u_j,-u_j)$ of critical points with $I(u_j)\to \infty$ as $j\to \infty$ by using the Symmetric mountain pass theorem \cite{Struwe}. We first need the following lemma:
\begin{lemma}\label{lems1}
Suppose that conditions
\eqref{H1} and \eqref{H2}  are satisfied. Then for any finite-dimensional subspace $\widetilde{X}\subset X$,
$
J_\lambda(u)\to-\infty, ~~\|u\|\to+\infty, ~u\in\widetilde{X}.
$
\end{lemma}
{\bf  Proof.}  Arguing indirectly, assume that there exists a sequence $(u_n)\subset \widetilde{X}$ such that
\begin{equation}\label{s1}
  \|u_n\|\to+\infty,~n\to+\infty~\text{and } J_{\lambda}(u_n)\geq B, \ \hbox{for all} \  n\in \mathbb{N},
\end{equation}
where $B\in \mathbb{R}$ is a fixed constant not depending on $n\in \mathbb{N}$.
Let $v_n=\frac{u_n}{\|u_n\|}$. Then
 it is obvious that $\|v_n\|=1$. Since dim $\widetilde{X}<+\infty$, there exists $v\in \widetilde{X}\backslash\{0\}$ such
that up to a subsequence,
$\|v_n-v\|\to 0,
\
\hbox{and}
\
v_n(x)\to v(x)
\quad
\hbox{a.e.}
\
x\in \Omega,
\
\hbox{as}
\
n\to +\infty.$
If $v(x)\neq0$, then $|u_n(x)|\to +\infty$ as $n\to +\infty$.
Clearly,  condition (\ref{H2}) implies condition
\begin{equation}\label{all1}
  \lim\limits_{|t|\to+\infty}\frac{F(x,t)}{|t|^{\frac{\gamma p^+}{r+1}}}=+\infty, ~\text{uniformly a.e. } x\in \Omega.
\end{equation}
By virtue of \eqref{all1},
{\small$$
\lim\limits_{n\to +\infty}\frac{F(x,u_n(x))}{\|u_n\|^{\frac{\gamma p^+}{r+1}}}=\lim\limits_{n\to +\infty}
\frac{F(x,u_n(x))}{|u_n|^{\frac{\gamma p^+}{r+1}}}|v_n|^{\frac{\gamma p^+}{r+1}}=+\infty,
x\in \Omega_0=\{x\in \Omega: v(x)\neq0\}.$$}Moreover, we can find $t_0>0$, such that
\begin{equation}\label{s2}
  \frac{F(x,t)}{|t|^{\frac{\gamma p^+}{r+1}}}\geq c>0,~~\hbox{for all} \  x\in \Omega ~\text{and } |t|>t_0.
\end{equation}
On the other hand, condition (\ref{H1}) implies that there exists a positive constant $C_1$ such that
\begin{equation}\label{s3}
  |F(x,t)|\leq C_1,\,\hbox{for all} \  (x,t)\in\Omega\times[-t_0,t_0].
\end{equation}
Then, by \eqref{s2} and \eqref{s3}, we deduce that there exists a constant $C_2\in\mathbb{R}$ such that
$ 
  F(x,t)\geq C_2,\,\hbox{for all} \  (x,t)\in\Omega\times\mathbb{R}.
$ 
From this, we conclude that
$ 
  {(F(x,u_{n}(x))-C_2)}{\|u_n\|^{-\frac{\gamma p^+}{r+1}}}\geq0,\,\hbox{for all} \  x\in\Omega \ \hbox{and} \  n\in\N,
$ 
which implies that
\begin{equation}\label{s4}
  \frac{F(x,u_{n}(x))}{|u_{n}(x)|^{\frac{\gamma p^+}{r+1}}}|v_{n}(x)|^{\frac{\gamma p^+}{r+1}}-\frac{C_2}{\|u_{n}\|^{\frac{\gamma p^+}{r+1}}}\geq0, \,\hbox{for all} \  x\in\Omega \ \hbox{and} \  n\in\N.
\end{equation}
Therefore using \eqref{s1} and \eqref{s4}, we have
\begin{eqnarray*}
0 &\leq& \lim_{n\to +\infty} J_{\lambda}(u_{n}(x)) \\
   &\leq&  \left\{
            \begin{array}{ll}
             \displaystyle\frac{\widehat{M}(1)}{(p^{+})^\gamma}-
             \displaystyle
             \lim_{n\to+\infty}
             \left[
              \frac{  \lambda \int_{\Omega}   |u_{n}(x)|^{p(x)}dx } {p^{-}\|u_n\|^{\frac{\gamma p^+}{r+1}}}
             +
              \frac{\left(
             \int_{\Omega} \frac{F(x,u_{n}(x))}{\|u_n\|^{\frac{\gamma p^+}{r+1}}}dx \right)^{r+1}}{r+1}
             \right],
              \hbox{if} ~\lambda\leq 0 \\
               \displaystyle\frac{\widehat{M}(1)}{(p^{+})^\gamma}-
               \displaystyle\lim_{n\to +\infty}\frac{1}{r+1}\left(\int_{\Omega}\frac{F(x,u_{n}(x))}{\|u_n\|^{\frac{\gamma p^+}{r+1}}}dx \right)^{r+1},
               \hbox{if}~ \lambda>0
            \end{array}
          \right.\\
    &\leq&  \left\{
            \begin{array}{ll}
             \displaystyle\frac{\widehat{M}(1)}{(p^{+})^\gamma}-
             \displaystyle\lim_{n\to +\infty}\left[\frac{\lambda C_{3}\|u_{n}\|^{p^+}}{p^{-}\|u_n\|^{\frac{\gamma p^+}{r+1}}}
             +
             \frac{\left(\int_{\Omega}\frac{F(x,u_{n}(x))-C_2}{\|u_n\|^{\frac{\gamma p^+}{r+1}}}dx \right)^{r+1}}{r+1}
             \right],
             \hbox{if} ~\lambda\leq 0 \\
               \displaystyle\frac{\widehat{M}(1)}{(p^{+})^\gamma}-
               \displaystyle\lim_{n\to +\infty}\frac{1}{r+1}\left(\int_{\Omega}\frac{F(x,u_{n}(x))-C_2}{\|u_n\|^{\frac{\gamma p^+}{r+1}}}dx \right)^{r+1},
               \hbox{if}~ \lambda>0
            \end{array}
          \right.\\
    &\leq& \displaystyle\frac{\widehat{M}(1)}{(p^{+})^\gamma}-
    \displaystyle\lim_{n\to +\infty}\frac{1}{r+1}\left(\int_{\Omega}\frac{F(x,u_{n}(x))-C_2}{\|u_n\|^{\frac{\gamma p^+}{r+1}}}dx \right)^{r+1}\\
    &\leq& \displaystyle\frac{\widehat{M}(1)}{(p^{+})^\gamma}-
    \displaystyle\lim_{n\to +\infty}\frac{1}{r+1}\left(\int_{\Omega}\frac{F(x,u_{n}(x))-C_2}{|u_n(x)|^{\frac{\gamma p^+}{r+1}}}\,|v_{n}(x)|^{\frac{\gamma p^+}{r+1}}dx \right)^{r+1}
    \to  -\infty,
\end{eqnarray*}
which is a contradiction. The proof of Lemma \ref{lems1} is thus complete.
\qed   
\subsection*{\bf Proof of Theorem \ref{theorem.1.2}.}
Clearly, by condition \eqref{H4}, $J_\lambda$ is an even functional. Since $J_\lambda(0) = 0$, thanks to Lemmas \ref{PS}, \ref{geometryMPT},
\ref{lems1} and the Symmetric mountain pass theorem \cite{Struwe}, we deduce the existence of an unbounded sequence of weak solutions to problem \eqref{Prob1}.\qed
\section{Proof of Theorem \ref{theorem.1.3}}\label{sec6}
In this part, we will prove Theorem \ref{theorem.1.3} by using  Clarke's theorem \cite{DC} which will be stated below. To this end, let us begin by defining the notion of genus and its basic properties.

 Let $\Sigma$ be the class of closed subset $A$ of $X\setminus 0$ such that $A=-A,$ i.e. symmetric with respect to the origin. Recall that for $A\in \Sigma$, the genus $\gamma(A)$ is defined as the least integer $k$ such that there exists an odd function $f\in C(X,\mathbb{R}^{k}\setminus 0).$ Moreover, if such function  does not exist then $\gamma(A)=\infty$ and by convenience, $\gamma(\emptyset)=0.$
 
It's well known that in general, the computation of the genus is a difficult task. Often, it suffices to use some estimates which can be given by comparison with sets whose genus is known as for example the sphere. We shall use the definition of the genus from \cite{MAK}.

Consider now  $\Sigma_{k}=\{A\in \Sigma, \gamma(A)\geq k\},
\
k\in \mathbb{N},$
and
$c_{k}:=\inf_{A\in \Sigma_{k}}\sup_{u\in A}I(u).$
We have $-\infty<c_{1}\leq c_{2}\leq \cdots \leq c_{k+1}\leq \cdots.$
Moreover, in order to prove Theorem \ref{theorem.1.3}, we use  Clarke's theorem.
\begin{theorem}[\cite{DC}]\label{Clarke's theorem}
Let $J\in C^{1}(X,\mathbb{R})$ be a functional satisfying the following conditions
\begin{itemize}
\item $(i)$  $J$ satisfies the $(PS)$ condition.
\item $(ii)$ $J$ is bounded from below and even.
\item $(iii)$ There exists a compact set $K\in \mathcal{A}$ such that $\gamma(A)=k$ and $\sup_{x\in K}J(x)<J(0).$
\end{itemize}Then $J$ possesses at least $k$ pairs of distinct critical points, and their corresponding critical values $c_{k}<0$ such that $\lim_{k\to \infty }c_{k}=0$ are less than $J(0).$
\end{theorem}
In order to get the infinity of solutions, we shall use Theorem \ref{Clarke's theorem}.
Since $X$ is a separable reflexive Banach space, there exist $(e_n)\subset X$ and $(e^{\star}_{n})\subset X^{\star}$ such that
$$\langle e^{\star}_{n}, e_m \rangle=\delta_{nm}=\begin{cases}
1~\hbox{if}~n=m\\
0~\hbox{if}~n\neq m,
\end{cases}
\quad
X=\overline{\hbox{span}\{e_n,~ n=1,2,\cdots,\}}, ~ X^{\star}=\overline{\hbox{span}\{e^{\star}_{n},~ n=1,2,\cdots,\}}.$$
For each $k\in \mathbb{N}$, consider the subspace
$X_{k}=\hbox{span}\{e_1,\cdots, e_k\} \subset
X=W_{0}^{1,p(x)}(\Omega),
\
\hbox{ spanned by}
\
e_1, \cdots,e_k.$ It is well known that $X_k\hookrightarrow L^{\delta(x)}(\Omega),$
 continuously for
 $1<\delta(x)<p^{\star}.$ Moreover, the norms in $X$ and $L^{\delta(x)}(\Omega)$ are equivalent in $X_k$.
Furthermore, by using  condition (\ref{H2}), we have $|F(x,u)|\geq C_{1}|u|^{\theta}-C_{2},$ hence we get
{\small$$J_{\lambda}(u)\leq \frac{\widehat{M}(1)}{{p^{-}}^{\gamma}}\left[\int_{\Omega}|\nabla u|^{p(x)}dx\right]^{\gamma} -\frac{\lambda}{p^{+}}\int_{\Omega}|u|^{p(x)}dx-\frac{C_{3}}{r+1}\left[\int_{\Omega}|u|^{(r+1)\theta}dx\right]+\frac{C_4 |\Omega|^{r+1}}{r+1}.$$}If $\|u\|$ is small enough, then we have $$\int_{\Omega}|\nabla u|^{p(x)}dx\leq \|u\|^{p^-}
\
\hbox{ and}
\
-|u|^{p^+}_{p(x)}\geq -\int_{\Omega}|u|^{p(x)}dx.$$ By using the equivalence of the norms in $X_k$, we deduce that
$$-C(k)\|u\|^{p^{+}}\geq -\int_{\Omega}|u|^{p(x)}dx,$$
where $C(k)$ is a positive constant.
Consequently, we get
$$J_{\lambda}(u)\leq \frac{\widehat{M}(1)}{{p^{-}}^{\gamma}}\|u\|^{\gamma p^{-}}-\frac{\lambda C(k)}{p^{+}}\|u\|^{p^+}-\tilde{C}(k)\|u\|^{(r+1)\theta}+C_{5}.$$
Hence, we have
{\small$$J_{\lambda}(u)
\leq
\|u\|^{(r+1)\theta}
\left[
\frac{\widehat{M}(1) \|u\|^{\gamma p^{-}-(r+1)\theta}}{{p^{-}}^{\gamma}}
 -
\frac{\lambda C(k) \|u\|^{p^{+}-(r+1)\theta}}{p^{+}}
+
\frac{C_{5}}{\|u\|^{(r+1)\theta}}
-
\tilde{C}(k)
\right].$$}Let $R$ be a positive constant such that $$\frac{\widehat{M}(1)}{{p^{-}}^{\gamma}}\|u\|^{\gamma p^{-}-(r+1)\theta} -\frac{\lambda C(k)}{p^{+}}\|u\|^{p^{+}-(r+1)\theta}+C_{5}\|u\|^{-(r+1)\theta}
\leq \tilde{C}(k).$$
Let  $0<r_{0}<R$ and consider the set $K=\{u\in X_{k};~ \|u\|=r_{0}\}$. Then
{\small\begin{eqnarray*}
J_{\lambda}(u)
&\leq&
r_{0}^{(r+1)\theta}\left[\frac{\widehat{M}(1)}{{p^{-}}^{\gamma}}r_{0}^{\gamma p^{-}-(r+1)\theta} -
\frac{\lambda C(k)}{p^{+}}r_{0}^{p^{+}-(r+1)\theta}+
C_{5}r_{0}^{-(r+1)\theta}-
\tilde{C}(k) \right]\\
&\leq&
R^{(r+1)\theta}
\left[
\frac{\widehat{M}(1) R^{\gamma p^{-}-(r+1)\theta}}{{p^{-}}^{\gamma}}
 -
 \frac{\lambda C(k) |g|_{\infty} R^{p^{+}-(r+1)\theta}}{p^{+}}
 +
 \frac{C_{5}}{R^{(r+1)\theta}}-\tilde{C}(k) \right]
< 0=J_{\lambda}(0),
\end{eqnarray*}}which implies that $\sup_{K}J_{\lambda}(u)<0=J_{\lambda}(0).$
Since $X_{k}$ and $\mathbb{R}^{k}$ are isomorphic and $K$ and $S^{k-1}$ are homomorphic, it follows that $\gamma(K)=k$. The Clarke theorem \ref{Clarke's theorem} shows that  problem \eqref{Prob1} admits at least $k$ pairs of distinct critical points, and their corresponding critical values $c_{k}<0$ such that $\lim_{k\to \infty }c_{k}=0$ are less than $J_{\lambda}(0).$ If $k$ is chosen arbitrary then problem \eqref{Prob1} possesses infinitely many critical points.
\begin{lemma}
For each $n\in \mathbb{N}$, there exists $\varepsilon >0$ such that
$\gamma(A^{-\varepsilon}_{\lambda})\geq n,
\
\hbox{where}
\
A^{-\varepsilon}_{\lambda}=\{u\in X;~ J_{\lambda}(u)\leq -\varepsilon\}.$
\end{lemma}
{\bf  Proof.} 
Consider $X_{n}$ be a subspace of $X$ of dimension $n$ and any $u\in X_{n}$ such that $\|u\|=1$ and $0<t<R$. Then we have
\begin{eqnarray*}
J_{\lambda}(tu)&\leq& \frac{\widehat{M}(1)t^{\gamma p^{-}}}{{p^{-}}^{\gamma}}\|u\|^{\gamma p^{-}}-\frac{\lambda C(k)t^{p^+}}{p^{+}}\|u\|^{p^+}-\tilde{C}(k)t^{(r+1)\theta}\|u\|^{(r+1)\theta}+C_{5}\\
&\leq& \frac{\widehat{M}(1)t^{\gamma p^{-}}}{{p^{-}}^{\gamma}}-\frac{\lambda C(k)t^{p^+}}{p^{+}}-\tilde{C}(k)t^{(r+1)\theta}+C_{5}.
\end{eqnarray*}
If $J_{\lambda}(tu)\to -\infty, \
\gamma p^{-}\leq \gamma p^{+}<(r+1)\theta , \ 0<t<R,$ then there exist $t_{0}>0$ and $\varepsilon >0$ such that
$J_{\lambda}(t_{0}u)<-\varepsilon,~ u\in X_{n},~ \|u\|=1.$
Consider now the sphere $S_{t_{0},n}=\{u\in X_n,~ \|u\|=t_{0}\}.$ Then $S_{t_{0},n}\subset A_{\lambda}^{-\varepsilon}$ and by the properties of the genus,  $\gamma(A_{\lambda}^{-\varepsilon})\geq \gamma{S_{t_{0},n}}=n$.
\qed   
\begin{lemma}
Let $\Sigma=\{A\subset X\setminus\{0\}
\mid
~A \ \hbox{is closed and} \ A=-A\},
$
and
$
\Sigma_{k}=\{A\in \Sigma
\mid
~ \gamma(A)\geq k\}.$ Then
$c_{k}=\inf_{A\in \Sigma_{k}}\sup_{u\in A}J_{\lambda}(u)$ is a negative critical value of $J_{\lambda}$ and if $c=c_{k}=\cdots=c_{k+r},$ then $\gamma(K_c)\geq r+1,$
where
$K_{c}=\{u\in X;~ J_{\lambda}(u)=c;~J'_{\lambda}(u)=0\}.$
\end{lemma}
{\bf  Proof.} 
First, we claim that $-\infty <c_{k}<\infty$. By the previous lemma, we know that for each $k\in \mathbb{N}$, there exists $\varepsilon >0$ such that $\gamma(A_{\lambda}^{-\varepsilon})\geq k$ or $A_{\lambda}^{-\varepsilon}\in \Sigma_{k}.$ Then
either
$c_{k}\leq \sup_{u\in A_{\lambda}^{-\varepsilon}}J_{\lambda}(u)\leq -\varepsilon(k)<0,\
\hbox{ for all}
\
K,$ or
$J_{\lambda}\
\hbox{ is bounded from below, hence}
\
c_{k}>-\infty,\
\hbox{ for all}
\
k\in \mathbb{N}.$
Since $c<0$, $J_{\lambda}$ satisfies the (PS) condition at level $c$, $K_c$ is compact and symmetric, it follows that $\gamma(K_c)$ is well-defined.

Let us now assume that $c=c_{k}=c_{k+1}=\cdots=c_{k+r}$
 and
 $\gamma(K_c)<r+1.$
By the properties of the genus, there exists a neighborhood $K$ of $K_c$ such that $\gamma(K)=\gamma(K_c)<r+1$. Moreover, based on the Deformation lemma \cite{PHR}, there exists an odd homomorphism
$\widehat{\eta}:X\rightarrow X$
such that
$\widehat{\eta}(A_{\lambda}^{c+\beta}\setminus K)\subset A_{\lambda}^{c-\beta},$
where
$ 0<\beta<-c.$ 

The functional $J_{\lambda}$ satisfies the (PS) in $A_{\lambda}^{0}$. Furthermore, by definition, we have
$c=c_{k+r}=\inf_{a\in \Sigma_{k+r}}\sup_{u\in A}J_{\lambda}(u).$
Then there exists $A\in \Sigma_{k+r}$ such that $\sup_{u\in A}J_{\lambda}(u)<c+\beta,$ which means that $A\subset A_{\lambda}^{c+\beta}$
 and
 $ \widehat{\eta}(A\setminus K)\subset \widehat{\eta}(A_{\lambda}^{c+\beta}\setminus K)\subset A_{\lambda}^{c-\beta}.$
Hence, we conclude that
$\gamma(\widehat{\eta}(\overline{A\setminus k}))\geq \gamma(\overline{A\setminus K})\geq \gamma(A)-\gamma(K)\geq (k+r)-r=k,$
i.e.,
$\widehat{\eta}(\overline{A\setminus K})\in \Sigma_{k},$
 hence
 $\sup_{u\in \widehat{\eta}(\overline{A\setminus K}) }J_{\lambda}\geq c_{k}=c,$ which yields a contradiction. Therefore, if $c=c_{k}=\cdots=c_{k+r},$ then $\gamma(K_c)\geq r+1$.
\qed   
\begin{remark}
We note that if $c_k$ is a critical value, then $\gamma(K_{c_{k}})\geq 1$ and $K_{c_{k}}$ is nonempty for all $k\in \mathbb{N}.$ In addition, if the points $c_{k}$ are not all distinct, then $\gamma(K_{c})> 1,$  $K_c$ is an infinite subspace, and  problem \eqref{Prob1} admits infinitely many critical points with negative energy.
\end{remark}

\section*{Acknowledgments}
Hamdani was supported by the Tunisian Military Research Center for Science and Technology Laboratory LR19DN01.
He also expresses his deepest gratitude to the Military Aeronautical Specialities School, Sfax (ESA) for providing an excellent
atmosphere for  work. Repov\v{s} was supported by the Slovenian Research Agency grants P1-0292, N1-0114 and N1-0083.
The authors wish to acknowledge the referees for several useful comments and valuable suggestions which have helped  improve the presentation.


\begin{thebibliography}{999}
\bibliographystyle{alpha}
\bibitem{A}
M. Allaoui, 
\emph{Existence results for a class of $p(x)$-Kirchhoff problems}, 
Studia Sci. Math. Hungar., $\mathbf{54}$(3), (2017) 316-331.

\bibitem{AD}
M. Allaoui  and A. Ourraoui,
\emph{Existence results for a class of $p(x)$-Kirchhoff problem with a singular weight}, 
Mediterr. J. Math., $\mathbf{13}$(2), (2016) 677-686.

\bibitem{AB}
C.O. Alves and  T. Boudjeriou, 
\emph{Existence of solution for a class of heat equation involving the p(x)-Laplacian with triple regime}, 
ZAMP, $\mathbf{72}$(1), (2021) 1-18.

\bibitem{AR} 
A. Ambrosetti and P.H. Rabinowitz, 
\emph{Dual variational methods in critical point theory and applications},
J. Funct. Anal., $\mathbf{14}$, (1973)  349-381.

\bibitem{BRW} 
A. Bahrouni, V.D. R\u{a}dulescu and P. Winkert, 
\emph{A critical point theorem for perturbed functionals and low perturbations of differential and nonlocal systems}, 
Adv. Nonlinear Stud.,  $\mathbf{20}$(3), (2020) 663-674.

\bibitem{B1}
C.J. Batkam,
\emph{An elliptic equation under the effect of two nonlocal terms}, 
Math. Methods Appl. Sci., $\mathbf{39}$(6), (2016) 1535-1547.

\bibitem{B2}
C.J. Batkam, 
\emph{Multiple sign-changing solutions to a class of Kirchhoff type problems},	
arXiv:1501.05733 [math.AP]

\bibitem{BBE}
Z. Binlin, G. Molica Bisci and R. Servadei,,
\emph{Superlinear nonlocal fractional problems with infinitely many solutions}, 
Nonlinearity, $\mathbf{28}$(7), (2015)  2247.

\bibitem{Chabrowski}
J. Chabrowski, 
\emph{On bi-nonlocal problem for elliptic equations with Neumann boundary conditions}, 
J. Anal. Math., $\mathbf{134}$, (2018) 303-334.

\bibitem{Choudhuri}
D. Choudhuri,
\emph{ Existence and H\"{o}lder regularity of infinitely many solutions to a p-Kirchhoff type problem
involving a singular nonlinearity without the Ambrosetti-Rabinowitz (AR) condition},
ZAMP,  $\mathbf{72}$(36), (2021).

\bibitem{DC}
D.C. Clarke, 
\emph{A variant of the Lusternik-Schnirelman theory,}
Indiana Univ.Math.J.,  $\mathbf{22}$, (1972) 65-74.

\bibitem{CC1}
F.J.S. Corr\^{e}a and   A.C.D.R. Costa, 
\emph{A variational approach for a bi-nonlocal elliptic problem involving the p(x)-Laplacian and non-linearity with non-standard growth}, Glasgow Math. J., $\mathbf{56}$(2), (2014) 317-333.

\bibitem{CC2}
F.J.S. Corr\^{e}a and   A.C.D.R. Costa, 
\emph{On a bi-non-local $p(x)$-Kirchhoff equation via Krasnoselskiis genus}, 
Math. Meth. Appl. Sci., $\mathbf{38}$(1), (2014) 87-93.

\bibitem{CF}
F.J.S. Corr\^{e}a and  G.M. Figueiredo, 
\emph{Existence and multiplicity of nontrivial solutions for a bi-nonlocal equation}, 
Adv. Differential Equ., {\bf 18}(5/6), (2013) 587-608.

\bibitem{FZ2}
X.-L. Fan and Q.-H. Zhang,
\emph{Existence of solutions for p(x)-Laplacian Dirichlet problem},
Nonlinear Anal., $\mathbf{52}$(8), (2003) 1843-1852.

\bibitem{H2}
M.K. Hamdani, 
\emph{On a nonlocal asymmetric Kirchhoff problem,} 
Asian-European J. Math., $\mathbf{13}$(5), (2019) 2030001.

\bibitem{HHMR}
M.K. Hamdani, A. Harrabi, F. Mtiri and D.D. Repov\v{s}, 
\emph{Existence and multiplicity results for a new p(x)-Kirchhoff problem},  
Nonlinear Anal., $\mathbf{190},$ (2020) 111598.

\bibitem{HR}
M.K. Hamdani and D.D. Repov\v{s}, 
\emph{Existence of solutions for systems arising in electromagnetism},
J. Math. Anal. Appl., $\mathbf{486}$(2),  (2020) 123898.

\bibitem{J}
Y. Jalilian, 
\emph{Infinitely many solutions for a bi-nonlocal equation with sign changing weight functions},
Bull. Iranian Math. Soc., $\mathbf{42}$(3), (2016) 611-626.

\bibitem{K1883}
G. Kirchhoff,
\emph{Mechanik}, 
Teubner, Leipzig, 1883.

\bibitem{MAK}
M. A. Krasnoselskii, 
\emph{Topological Methods in the Theory of Nonlinear Integral Equations,}
MacMillan, New York, (1964).

\bibitem{Lions}
J.L. Lions,
\emph{On some questions in boundary value problems of mathematical physics},
North-Holland Math. Stud., $\mathbf{30}$, (1978)
284-346

\bibitem{PRD}
N.S. Papageorgiou, V.D. R\u{a}dulescu and D.D. Repov\v{s},
\emph{Nonlinear Analysis - Theory and Methods},
Springer Monographs in Mathematics, Springer, Cham, 2019.

\bibitem{PHR}
P.H. Rabinowitz,
\emph{Minimax methods in critical point theory with applications to differential equations,}
CBMS Reg. Conf. Ser. Math.,  $\mathbf{65}$ (1984).

\bibitem{RD}
V.D. R\u{a}dulescu and D.D. Repov\v{s},
\emph{Partial Differential Equations with Variable Exponents: Variational Methods and Qualitative Analysis},
 CRC Press, Boca Raton, 2015.

\bibitem{Struwe}  
M. Struwe,
\emph{Variational Methods, Applications to Nonlinear Partial Differential Equations and Hamiltonian Systems}, Second ed., 
Springer-Verlag, Berlin, 1996.

\bibitem{Willem}
M. Willem, 
\emph{Minimax Theorems}, 
Progress in Nonlinear Differential Equations and their Applications, Vol.~24,
Birkh\"{a}user Boston, Inc., Boston, 1996.

\bibitem{Y}
Z. Yucedag, 
\emph{Existence of solutions for p(x)-Laplacian equations without Ambrosetti-Rabinowitz type condition}, 
Bull. Malays. Math. Sci. Soc., $\mathbf{38}$(3), (2015) 1023-1033.

\bibitem{ZGC}
B. L. Zhang, B. Ge and X.-F. Cao, 
\emph{Multiple Solutions for a Class of New p(x)-Kirchhoff Problem without the Ambrosetti-Rabinowitz conditions}, 
Mathematics, $\mathbf{8}$(11), (2020) 2068.

\bibitem{ZFB}
J. Zuo, A. Fiscella and A. Bahrouni, 
\emph{Existence and multiplicity results for $p(\cdot)\&q(\cdot)$ fractional Choquard problems with variable order}, 
Complex Var. Elliptic Equ., $\mathbf{67}$(2), (2020) 500-516.
\end{thebibliography}
\end{document}